# Nonexistence of simple hyperbolic blow-up for the quasi-geostrophic equation

By Diego Cordoba

*In memory of Maria del Carmen Gazolaz*

## 1. Introduction

The problem we are concerned with is whether singularities form in finite time in incompressible fluid flows. It is well known that the answer is "no" in the case of Euler and Navier-Stokes equations in dimension two. In dimension three it is still an open problem for these equations.

In this paper we focus on a two-dimensional active scalar model for the 3D Euler vorticity equation. Constantin, Majda and Tabak [7] suggested, by studying rigorous theorems and detailed numerical experiments, a general principle: "If the level set topology in the temperature field for the 2D quasi-geostrophic active scalar in the region of strong scalar gradients does not contain a hyperbolic saddle, then no finite time singularity is possible."

Numerical simulations showed evidence of singular behavior when the geometry of the level sets of the active scalar contain a hyperbolic saddle. There is a naturally associated notion of simple hyperbolic saddle breakdown. The main theorem we present in this paper shows that such breakdown cannot occur in finite time. We also show that the angle of the saddle cannot close in finite time and it cannot be faster than a double exponential in time. Using the same techniques, we see that analogous results hold for incompressible 2D and 3D Euler.

These results were announced in [9], but with a slight difference in the definition of a simple hyperbolic saddle. The definition given in Section 4 generalizes the one given in the announcement. See also Constantin [4], discussed in Section 7, Remark 5 below.

The whole work described in this paper is basically part of the author's thesis. I am particularly grateful to my thesis advisor Charles Fefferman for his attention, support, guidance and advice. I am indebted to D. Christodoulou and P. Constantin for helpful corrections and suggestions. I wish to thank A. Majda for suggesting the subject and E. Tabak for discussions and com-ments on the first version of the paper. The author was supported by the Sloan Foundation during part of the writing of this paper.

## 2. Quasi-geostrophic equation as a 2D model for 3D Euler

The vorticity equation for incompressible 3D Euler is the following

$$(\partial_t + u \cdot \nabla)\omega = (\nabla u) \cdot \omega, \tag{1}$$

where $\omega$ is the vorticity and the velocity is divergence free. Using the Biot-Savart law we can recover the velocity from the vorticity by

$$u(x,t) = \frac{1}{4\pi} \int_{R^3} \frac{y \times \omega(x+y,t)}{|y|^3} dy.$$

Beale, Kato and Majda [1] showed that a necessary and sufficient condition to have a singularity at time $T$ is that

$$\int_0^T \|\omega(t)\|_\infty dt = +\infty.$$

Based on this result, Constantin, Fefferman and Majda [5] proved that if the direction field $\xi(x) = \frac{\omega(x)}{|\omega(x)|}$ is smooth in and near regions of high $|\omega|$ then blow-up does not occur.

Another way to understand the problem is by constructing and studying models in lower dimensions. Constantin, Lax and Majda [6] studied a one-dimensional model for the 3D Euler vorticity equation,

$$\frac{\partial \omega}{\partial t} = H(\omega)\omega.$$

Here $H$ is the Hilbert transform and the velocity is defined by $u = \int_{-\infty}^x \omega(y)dy$. They show that there are solutions with breakdown.

A 2D model of the quasi-geostrophic equation was studied in [7]. This equation represents a geophysical model [14], where the conjectured singularities describe the formation of sharp fronts between masses of hot and cold air.

The 2D quasi-geostrophic active scalar equations are

$$(\partial_t + u \cdot \nabla)\theta = 0 \tag{2}$$

$$u = \nabla^\perp \psi \quad \text{where} \quad \theta = -(-\triangle)^{\frac{1}{2}}\psi$$

and the initial condition $\theta(x,0) = \theta_0$.

Here $\theta = \theta(x,t)$ with $x \in \mathbf{R}^2, t \in \mathbf{R}^+$ is a scalar, $u$ is the velocity, and $\psi$ is the stream function.

In [7], [14], $\theta$ is interpreted as a potential temperature. I am grateful to R. de la Llave for helpful comments, and for pointing out to me that in the quasi-geostrophic equation, $\theta$ may also be regarded as a vorticity.

In [7] it was shown that the level sets of $\theta$ are analogous to the vortex lines for 3D Euler and that there is a geometric and analytic similarity between both equations. It was also shown there that results analogous to the ones mentioned above for 3D Euler also hold for the 2D active scalar equation:

(i) By differentiating the equation we obtain an equation similar to (1) but in this case, instead of the vorticity, there is the gradient perpendicular to $\theta$:
$$(\partial_t + u \cdot \nabla) \nabla^\perp \theta = (\nabla u) \cdot \nabla^\perp \theta.$$

The stream function can be written
$$\psi(y_1, y_2) = -\int_{R^2} \frac{\theta(x+y)}{|x|} dx;$$

thus
$$u(x,t) = -\int_{R^2} \frac{\nabla^\perp \theta(x+y,t)}{|y|} dy,$$

and obviously the velocity is divergence free.

(ii) If $\theta_o(x) \in H^k(\mathbf{R}^2)$ $k \geq 3$, then a necessary and sufficient condition for having a singularity at time $T$ is that

(3) $$\int_0^T \|\nabla^\perp \theta\|_\infty dt = +\infty.$$

(iii) If the direction field $\xi(x) = \frac{\nabla^\perp \theta}{|\nabla^\perp \theta|}$ is smooth in and near regions of high $|\nabla^\perp \theta|$ then blow-up does not occur.

Another similarity is that in both cases the kinetic energy is conserved for all time.

### 3. Numerical simulations

It is not known if the quasi-geostrophic equation develops singularities in finite time. There have been several numerical attempts to find a candidate for initial data such that a strong singular behavior develops in finite time. In [7] the initial data
$$\theta(x,0) = \sin(x_1)\sin(x_2) + \cos(x_2)$$
were studied numerically and found to develop a strong front. An empirical asymptotic fit for $\max|\nabla^\perp \theta|$ was obtained behaving like $(T_* - t)^{-1.66}$ with $T_* \cong 8.25$.

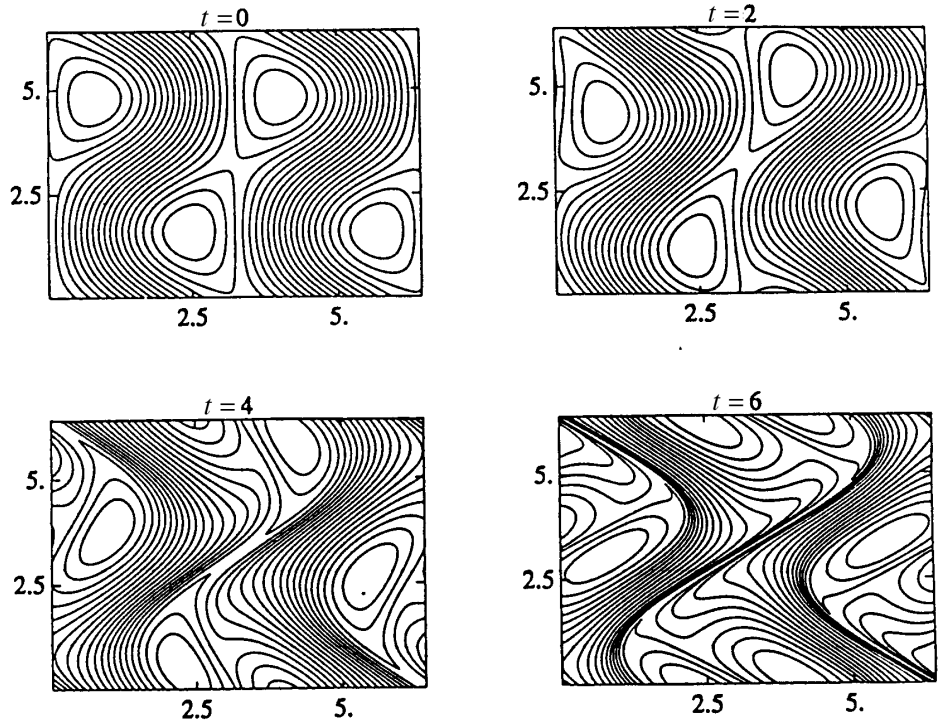

Figure 1. Numerical simulations

Figure 1 represents the evolution of the initial data mentioned above at times $t = 0, 2, 4, 6$. The lines are the level sets of $\theta$. Is clear that for $t = 0$ the level sets contain a hyperbolic saddle. As time evolves, Figure 1 shows that the saddle is closing very fast. Therefore there is a strong front formation. The author thanks P. Constantin, A. Majda and E. Tabak for permission to reproduce these pictures from their paper ([7]).

Ohkitani and Yamada [13] gave another interpretation of this particular front. They suggested that the maximum gradient does not go to infinity in finite time, but rather goes to infinity at a double exponential rate in time. Recently, Constantin, Nie, and Schorghofer [8] made careful measurements of stretching rates and collected substantial numerical evidence against a singularity for this particular case.

No other initial data that develop a front faster than the one mentioned above have been found.

## 4. Hyperbolic saddle scenario

To rigorously and analytically approach the problem where the level sets of the active scalar field contain a hyperbolic saddle, we first need a precise definition of a hyperbolic saddle.

In this section we define a naturally associated notion of simple hyperbolic blow-up scenario and we rule out a finite time singularity.

*Definition* 1. A simple hyperbolic saddle in a neighborhood $U$ of the origin is the set of curves $\rho = \text{const}$ where

$$\rho = (y_1\beta(t) + y_2)(y_1\delta(t) - y_2),$$

and there is a nonlinear time-dependent coordinate change

$$y_1 = F_1(x_1, x_2, t)$$

$$y_2 = F_2(x_1, x_2, t)$$

with $\beta(t), \delta(t) \in \mathcal{C}^1([0, T_*))$, $F_i \in \mathcal{C}^2(\overline{U} \times [0, T_*])$, $|\beta|, |\delta| \leq C$, $\beta(t) + \delta(t) \geq 0$, $|\det\frac{\partial F_i}{\partial x_j}| \geq c > 0$ whenever $x \in U$, $t \in [0, T_*]$.

*Remark* 1. The saddle is allowed to rotate and dilate with respect to time. The center of the saddle can move in $U$ with time.

*Remark* 2. The angle of opening of the saddle is $\gamma \simeq \beta + \delta$. The restrictions $\beta(t) + \delta(t) \geq 0$ and $|\beta|, |\delta| \leq C$ are unimportant. See also [11].

*Remark* 3. The definition given in [9] for $\rho$ was $\rho = y_1 y_2 - \cot\alpha \cdot y_2^2$ where only one branch, $y_2 = \tan\alpha \cdot y_1$, was allowed to close and the other one moved smoothly. In the present paper both branches ($y_2 = -\beta(t)y_1, y_2 = \delta(t)y_1$) are allowed to close at different rates. The results of [9] are therefore a particular case of the one presented here.

The possible singularity in this scenario is due to $\gamma(t)$ becoming zero in finite time. The following theorem will show that this is not possible and $\gamma(t)$ can go to zero at most as a double exponential.

THEOREM 1. *Let $\theta(x_1, x_2, t)$ be a smooth solution of (2) defined for $0 \leq t < T_*$, $(x_1, x_2) \in R^2$. Assume for $0 \leq t < T_*$ that $\theta$ is constant along the curves $\rho = \text{const}$ defined in Definition 1. Assume also, for each fixed $t$, that $\theta$ is not constant on any disc in $U$. Then $\lim_{t \to T_*} \gamma(t)$ exists and is not $0$.*

COROLLARY 1. *Let $\theta$ be as in Theorem 1, let $\xi = \frac{\nabla^\perp \theta}{|\nabla^\perp \theta|}$, and assume $|\nabla \xi| < C$ on $(R^2 \setminus U) \times [0, T_*]$. Then $\theta$ continues to some solution of (2) on $R^2 \times [0, T_* + \varepsilon]$ for some $\varepsilon > 0$.*

THEOREM 2. *Let $\theta(x_1, x_2, t)$, $\beta$, $\delta$, $U$ and $F_j$ be as in Theorem 1, but with $T_* = \infty$. Assume that the $C^2$ seminorms of $F_j$ are bounded for all time $t \in [0, \infty)$. Then*

$$|\log \log \frac{1}{\gamma(t)}| \leq (\text{constant}) \cdot t + \text{const}$$

*for all t.*

COROLLARY 2. *Let $\theta$ be as in Theorem 2. Let $\xi = \frac{\nabla^\perp \theta}{|\nabla^\perp \theta|}$. Assume $|\nabla \xi| \leq \Phi(t)$ on $(R^2 \setminus U) \times [0, \infty)$. Then*

$$|\nabla^\perp \theta| \leq e^{e^{c \int_0^t (e^{e^s} + \Phi(s))ds}}$$

*on $R^2 \times [0, \infty)$.*

*Remark* 4. In view of Theorem 2, the most natural example is $\Phi(t) = e^{e^t}$. Corollary 2 then shows that $|\nabla^\perp \theta|$ is bounded by a quadruple exponential.

## 5. Stream function under a change of variables

The purpose of this section is to obtain an expression for the stream function under a new set of variables $(\rho, \sigma)$. The variable $\rho = \rho(x_1, x_2, t)$ was defined in Definition 1, and $\sigma$ is defined as follows:

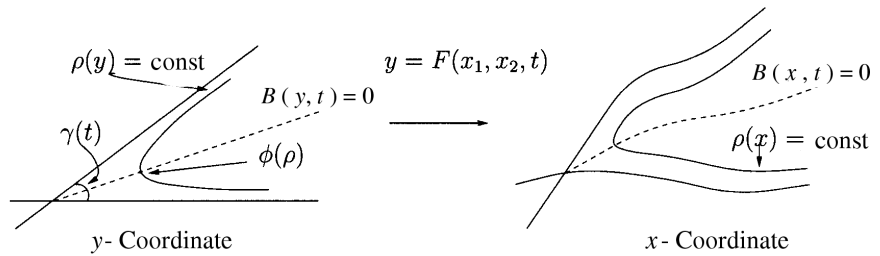

Figure 2. Change of variables

Let $B(y_1, y_2, t) = 0$ be the bisector of the angle $\gamma$, $\rho \geq 0$ and we denote $\phi(\rho)$ to be the intersection of $B = 0$ with $\rho$ as shown in Figure 2. We write $\tilde{\phi}(x)$ for $\phi(\rho(x))$.

We then define $\sigma(x)$ to be the solution of

(4) $$\exp(\sigma \nabla_x^\perp \rho)[\tilde{\phi}(x)] = x.$$

Note that, as $\sigma$ varies, $\exp(\sigma \nabla_x^\perp \rho)[\tilde{\phi}(x)]$ moves along the curve $\rho = \text{const}$. Thus $\sigma(x)$ is well defined for $x \in U$ and $\rho > 0$.

Let $\theta(x_1, x_2, t)$ be a smooth solution of (2) defined for $0 \le t < T_*$, $(x_1, x_2) \in R^2$. Assume for $0 \le t < T_*$ that $\theta$ is constant along the curves $\rho = $ const where $\rho = \rho(x_1, x_2, t)$. Hence we write $\theta(x, t) = \tilde{\theta}(\rho, t)$ for a function $\tilde{\theta}$. Also, $\psi(x,t) = \tilde{\psi}(\rho, \sigma, t)$ and $\nabla^\perp \theta$ can be expressed as

$$\nabla^\perp \theta = \frac{|\nabla^\perp \theta|}{r}(\zeta_1, \zeta_2),$$

where $r = \sqrt{\zeta_1^2 + \zeta_2^2}$ and

$$\zeta_1 = \frac{\partial x_1}{\partial \sigma} = -\frac{\partial \rho}{\partial x_2}$$

$$\zeta_2 = \frac{\partial x_2}{\partial \sigma} = \frac{\partial \rho}{\partial x_1}.$$

By making the change of variables in equation (2), we obtain

$$u \cdot \nabla_x \theta = \frac{\partial \tilde{\theta}}{\partial \rho}(u \cdot \nabla_x \rho) = -\frac{\partial \tilde{\theta}}{\partial \rho}\left(\frac{\partial \psi}{\partial x_2}\frac{\partial x_2}{\partial \sigma} + \frac{\partial \psi}{\partial x_1}\frac{\partial x_1}{\partial \sigma}\right);$$

thus

$$\frac{\partial \tilde{\theta}}{\partial t} + \frac{\partial \tilde{\theta}}{\partial \rho}\left(\frac{\partial \rho}{\partial t} - \frac{\partial \psi}{\partial \sigma}\right) = 0.$$

Taking into account that $\frac{\partial \tilde{\theta}}{\partial \rho}$ is independent of $\sigma$ and that $\theta$ is not constant in a disk, we see that the derivative with respect to time of $\rho$ along trajectories, i.e. $\frac{\partial \rho}{\partial t} - \frac{\partial \psi}{\partial \sigma}$, is equal to a function $(-H_1)$ that is also independent of $\sigma$. Therefore

(5) $$\frac{\partial \psi}{\partial \sigma} = \frac{\partial \rho}{\partial t} + H_1(\rho, t).$$

Furthermore, integrating with respect to $\sigma$ gives

(6) $$\psi(\rho, \sigma, t) = H_1(\rho, t) \cdot \sigma + \int_0^\sigma \frac{\partial \rho}{\partial t} d\sigma + H_2(\rho, t),$$

where $H_2(\rho, t) = \tilde{\psi}(\rho, 0, t)$.

## 6. Proofs

The main idea in the proof of both theorems is to estimate the difference of the value of the stream function at a point $p$ that lies in the branch of the saddle $y_2 = -\beta(t) y_1$ with the value of the stream function at $q$ that lies in the other branch $y_2 = \delta(t) y_1$. Both $p$ and $q$ have the same $y_1$ coordinate. We need two expressions of the stream function; one comes from the equality $\theta = -(-\triangle)^{\frac{1}{2}}\psi$ and the other one comes from the change of variables done in Section 5.

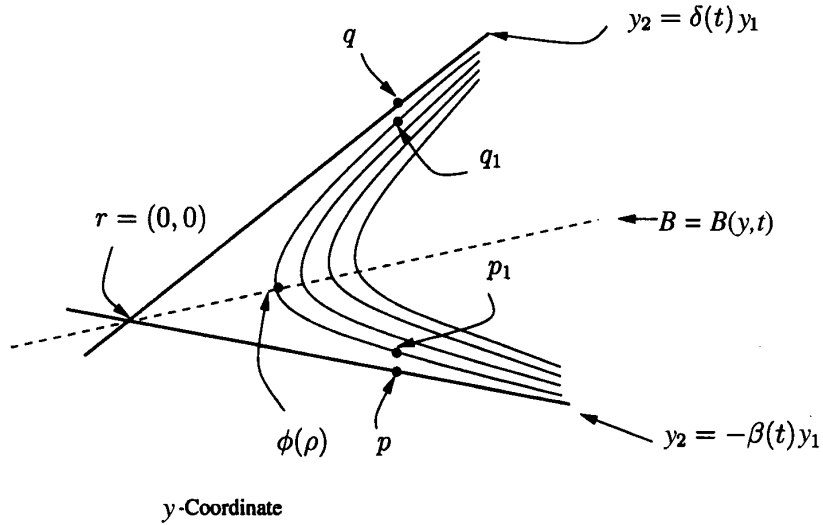

Figure 3. Defining $p$, $q$, $p_1$, $q_1$, and $r$

LEMMA 1. *Let $\theta$ be a solution of equation (2), $p$ and $q$ be as defined above and $\psi$ be given by $\psi = -(-\triangle)^{-\frac{1}{2}}\theta$. Then*

$$|\psi(p) - \psi(q)| \leq K_1|\gamma \cdot \log \gamma|$$

*where $K_1$ is a constant and $|p - q| \sim \gamma$.*

To prove this lemma we use the fact that $\theta$ is bounded for all time; this is because the derivative with respect to time along the Lagrangian trajectories is zero. We also recall that $\theta(x,t)$ belongs to $L^2$ for fixed $t$, with $L^2$ norm independent of $t$.

$$I = \psi(p) - \psi(q) = \int_{R^2} \theta(y)\left(\frac{1}{|y-p|} - \frac{1}{|y-q|}\right)dy.$$

Let us denote $\tau = |p - q|$. We split the integral $I$:

$$I(x) = \int_{|y-p|\leq 2\tau} + \int_{2\tau<|y-p|\leq k} + \int_{k<|y-p|} = I_1 + I_2 + I_3,$$

where $k$ is a fixed number.

We bound $I_1$ by

$$|I_1| \leq ||\theta||_{L^\infty} \cdot \int_{|y-p|\leq 2\tau} \left|\frac{1}{|y-p|} - \frac{1}{|y-q|}\right|dy$$

$$\leq C \cdot \int_{|y-p|\leq 2\tau} \left(\frac{1}{|y-p|} + \frac{1}{|y-q|}\right)dy$$

$$\leq C\tau.$$

We define $s$ to be a point in the line between $p$ and $q$. Then $|y-p| \leq 2|y-s|$, and $I_2$ can be estimated by

$$|I_2| \leq C\tau \cdot \int_{2\tau < |y-p| \leq k} \max_s |\nabla(\frac{1}{|y-s|})| dy$$
$$\leq C\tau \cdot \int_{2\tau < |y-p| \leq k} \max_s \frac{1}{|y-s|^2} dy$$
$$\leq C\tau \cdot |\log \tau|.$$

To estimate $I_3$ we recall that $\int_{R^2} |\theta|^2 dx$ is conserved for all time. It is easy to check that $|I_3| \leq C \cdot \tau$. Therefore $|I| \leq \tau \cdot |\log \tau|$.

To state and prove the following lemmas we have to define

$$\tilde{q}(y_1, t) = (y_1, \delta(t) y_1)$$
$$\tilde{p}(y_1, t) = (y_1, -\beta(t) y_1).$$

LEMMA 2. *With the same assumptions as in Theorem 1, let $\psi$ be given by expression* (6) *and $(p, q)$ defined as before. Then*

$$S_1 = \psi(q) - \psi(p) = \frac{d\delta}{dt} \cdot \int_0^{y_1} \frac{\tilde{y}_1}{D(\tilde{q}(\tilde{y}_1, t))} d\tilde{y}_1 + \frac{d\beta}{dt} \cdot \int_0^{y_1} \frac{\tilde{y}_1}{D(\tilde{p}(\tilde{y}_1, t))} d\tilde{y}_1 + O(\gamma)$$

*where $D = |\det \frac{\partial F_i}{\partial x_j}|$.*

*Proof.* We evaluate $\psi$ at the points $p_1 = (\rho, \sigma_1)$ and $q_1 = (\rho, \sigma_2)$ with $\sigma_1 \neq \sigma_2$, as in Figure 3

$$\psi(q_1) - \psi(p_1) = H_1(\rho, t) \cdot (\sigma_1 - \sigma_2) + \int_{\sigma_1}^{\sigma_2} \frac{\partial \rho}{\partial t} d\sigma.$$

The next step is to take the limits $p_1 \to p$ and $q_1 \to q$. That means $\rho \to 0$.

For the first term of the right-hand side of the equality we consider $t \in [0, T_*)$; hence the velocity is bounded as a function of $x$. By the construction of the set of variables $(\rho, \sigma)$ in Section 5, it follows that the $\sigma$-coordinates of $p$, $q$ grow at most as $\log \rho$ when $p_1 \to p$ and $q_1 \to q$. On the other hand $H_1$ is independent of $\sigma$ and is given by the following equation:

$$H_1(\rho, t) = -u_1 \frac{\partial \rho}{\partial x_1} - u_2 \frac{\partial \rho}{\partial x_2} - \frac{\partial \rho}{\partial t}.$$

For a fixed $t \in [0, T_*)$ we can bound $H_1$ by

$$|H_1(\rho, t)| \leq |u| \cdot |\nabla \rho| + |\frac{\partial \rho}{\partial t}|,$$

where

$$\frac{\partial \rho}{\partial x_i} = \frac{\partial \rho}{\partial y_1} \frac{\partial y_1}{\partial x_i} + \frac{\partial \rho}{\partial y_2} \frac{\partial y_2}{\partial x_i}$$

$$\frac{\partial \rho}{\partial t} = \frac{\partial \rho}{\partial y_1}\frac{\partial y_1}{\partial t} + \frac{\partial \rho}{\partial y_2}\frac{\partial y_2}{\partial t} + \frac{d(\delta - \beta)}{dt} y_1 \cdot y_2 + \frac{d(\beta\delta)}{dt} y_1^2.$$

We take into account the definition of $\rho = \rho(y_1, y_2, t)$, $y = F(x_1, x_2, t)$ and $F_i \in \mathcal{C}^2(\overline{U} \times [0, T_*])$; then

$$|\frac{\partial \rho}{\partial x_i}| \leq |y| \cdot (\text{const})$$

$$|\frac{\partial \rho}{\partial t}| \leq |y| \cdot (\text{const}).$$

Thus

$$|H_1(\rho, t)| \leq |y| \cdot (\text{const})$$

(these constants may depend on $t$), and, from the definition of $\rho$, we also know that $|y| \sim \rho^{\frac{1}{2}}$ when we approach the origin along the bisector $B$ with $\rho > 0$. Therefore $|H_1|$ is at most $\rho^{\frac{1}{2}}$ when $\rho \to 0$. This implies that

$$\lim_{\rho \to 0} H_1(\rho, t) \cdot (\sigma_1 - \sigma_2) = 0.$$

For the second term we obtain something different. We define

$$\Gamma = [(y_1, y_2) : \rho = \text{const}].$$

We have the following expressions for $\frac{\partial y_i}{\partial \sigma}$ by the change of variables

$$\frac{\partial y_i}{\partial \sigma} = \frac{\partial y_i}{\partial x_1}\frac{\partial x_1}{\partial \sigma} + \frac{\partial y_i}{\partial x_2}\frac{\partial x_2}{\partial \sigma} = -\frac{\partial y_i}{\partial x_1}\frac{\partial \rho}{\partial x_2} + \frac{\partial y_i}{\partial x_2}\frac{\partial \rho}{\partial x_1}$$
$$= -\frac{\partial y_i}{\partial x_1}\left(\frac{\partial \rho}{\partial y_1}\frac{\partial y_1}{\partial x_2} + \frac{\partial \rho}{\partial y_2}\frac{\partial y_2}{\partial x_2}\right) + \frac{\partial y_i}{\partial x_2}\left(\frac{\partial \rho}{\partial y_1}\frac{\partial y_1}{\partial x_1} + \frac{\partial \rho}{\partial y_2}\frac{\partial y_2}{\partial x_1}\right).$$

Then

$$\frac{\partial y_1}{\partial \sigma} = -D\frac{\partial \rho}{\partial y_2} \quad , \quad \frac{\partial y_2}{\partial \sigma} = D\frac{\partial \rho}{\partial y_1}$$
$$\frac{\partial y_1}{\partial \sigma}\frac{d\sigma}{dy_1} = 1 \quad , \quad \frac{\partial y_2}{\partial \sigma}\frac{d\sigma}{dy_2} = 1 \quad \text{on} \quad \Gamma.$$

Using these formulas for the integral on $\Gamma$ we find that

$$\int_{\sigma_1}^{\sigma_2} \frac{\partial \rho}{\partial t} d\sigma = \int_{\sigma_1}^{\sigma_2} \frac{\partial \rho}{\partial y_1}\frac{\partial y_1}{\partial t} d\sigma + \int_{\sigma_1}^{\sigma_2} \frac{\partial \rho}{\partial y_2}\frac{\partial y_2}{\partial t} d\sigma$$
$$+ \frac{d(\delta - \beta)}{dt}\int_{\sigma_1}^{\sigma_2} y_1 \cdot y_2 d\sigma + \frac{d(\beta\delta)}{dt}\int_{\sigma_1}^{\sigma_2} y_1^2 d\sigma$$
$$= [I_1(q) - I_1(p)] + \cdots + [I_4(q) - I_4(p)],$$

where

$$I_1 = -\int_0^{y_1} \frac{1}{D \cdot \frac{\partial \tilde{\rho}}{\partial \tilde{y}_2}} \cdot \frac{\partial \tilde{\rho}}{\partial \tilde{y}_1} \frac{\partial \tilde{y}_1}{\partial t} d\tilde{y}_1$$

$$I_2 = -\int_0^{y_1} \frac{1}{D} \frac{\partial \tilde{y}_2}{\partial t} d\tilde{y}_1$$

$$I_3 = -\frac{d(\delta - \beta)}{dt} \int_0^{y_1} \frac{\tilde{y}_1 \cdot \tilde{y}_2}{D \cdot \frac{\partial \tilde{\rho}}{\partial \tilde{y}_2}} d\tilde{y}_1$$

$$I_4 = -\frac{d(\beta\delta)}{dt} \int_0^{y_1} \frac{\tilde{y}_1^2}{D \cdot \frac{\partial \tilde{\rho}}{\partial \tilde{y}_2}} d\tilde{y}_1.$$

We obtain the following estimates for $I_i(q) - I_i(p)$:

*Case $I_1$.*

$$I_1(q) - I_1(p) = \int_0^{y_1} \frac{\delta}{D((\tilde{q}(\tilde{y}_1, t))} \frac{\partial \tilde{y}_1}{\partial t}(\tilde{q}(\tilde{y}_1, t)) d\tilde{y}_1$$
$$+ \int_0^{y_1} \frac{\beta}{D((\tilde{p}(\tilde{y}_1, t))} \frac{\partial \tilde{y}_1}{\partial t}(\tilde{p}(\tilde{y}_1, t)) d\tilde{y}_1$$
$$= (\delta + \beta) \int_0^{y_1} \frac{\frac{\partial \tilde{y}_1}{\partial t}(\tilde{q})}{D(\tilde{q})} d\tilde{y}_1 + \beta \int_0^{y_1} \left( \frac{\frac{\partial \tilde{y}_1}{\partial t}(\tilde{p})}{D(\tilde{p})} - \frac{\frac{\partial \tilde{y}_1}{\partial t}(\tilde{q})}{D(\tilde{q})} \right) d\tilde{y}_1$$
$$= O(\gamma).$$

*Case $I_2$.*

$$I_2(q) - I_2(p) = \int_0^{y_1} \left( \frac{\frac{\partial \tilde{y}_2}{\partial t}(\tilde{p}(\tilde{y}_1, t))}{D(\tilde{p}(\tilde{y}_1, t))} - \frac{\frac{\partial \tilde{y}_2}{\partial t}(\tilde{q}(\tilde{y}_1, t))}{D(\tilde{q}(\tilde{y}_1, t))} \right) d\tilde{y}_1$$
$$= O(\gamma).$$

*Case $I_3$.*

$$I_3^* = I_3(q) - I_3(p) = -\frac{d(\delta - \beta)}{dt} \cdot \frac{1}{\beta + \delta} \int_0^{y_1} \left( \frac{\beta}{D(\tilde{p}(\tilde{y}_1, t))} - \frac{\delta}{D(\tilde{q}(\tilde{y}_1, t))} \right) \tilde{y}_1 d\tilde{y}_1.$$

*Case $I_4$.*

$$I_4^* = I_4(q) - I_4(p) = \frac{d(\beta\delta)}{dt} \cdot \frac{1}{\beta + \delta} \int_0^{y_1} \left( \frac{1}{D(\tilde{q}(\tilde{y}_1, t))} + \frac{1}{D(\tilde{p}(\tilde{y}_1, t))} \right) \tilde{y}_1 d\tilde{y}_1.$$

Above we used the fact that

$$\frac{\partial \rho}{\partial y_2}(\tilde{q}(\tilde{y}_1, t)) = -\tilde{y}_1(\beta + \delta), \quad \frac{\partial \rho}{\partial y_2}(\tilde{p}(\tilde{y}_1, t)) = \tilde{y}_1(\beta + \delta)$$
$$\frac{\partial \rho}{\partial y_1}(\tilde{q}(\tilde{y}_1, t)) = \delta\tilde{y}_1(\beta + \delta), \quad \frac{\partial \rho}{\partial y_1}(\tilde{p}(\tilde{y}_1, t)) = \beta\tilde{y}_1(\beta + \delta),$$

and that the function $\frac{1}{D}\frac{\partial y_i}{\partial t}$ is $\mathcal{C}^1(\overline{U} \times [0, T_*])$.

Now we can complete the proof of Lemma 2 by adding $I_3^*$ to $I_4^*$

$$I_3^* + I_4^* = \frac{d\delta}{dt} \cdot \int_0^{y_1} \frac{\tilde{y}_1}{D(\tilde{q}(\tilde{y}_1, t))} d\tilde{y}_1 + \frac{d\beta}{dt} \cdot \int_0^{y_1} \frac{\tilde{y}_1}{D(\tilde{p}(\tilde{y}_1, t))} d\tilde{y}_1.$$

LEMMA 3. *Under the same assumptions as above and $(q, r)$ as in Figure 2, we then have*

$$S_2 = \psi(q) - \psi(r) = \frac{d\delta}{dt} \cdot \int_0^{y_1} \frac{\tilde{y}_1}{D(\tilde{q}(\tilde{y}_1, t))} d\tilde{y}_1 + E(x_1, x_2, t)$$

*where $E$ is bounded for all time.*

To prove the equality we again use the expression (6) for the stream function, and evaluate it at $\phi(\rho)$ and $q_1$, where these points belong to the same level set $\rho$. Then we take the limit $\rho \to 0$, which means that $\phi(\rho) \to r$ and $q_1 \to q$. The equality follows from using the same steps as in the proof of Lemma 2. The first two cases, $I_1$ and $I_2$, give us the function $E(x_1, x_2, t)$, and the other two cases give us the other term of the right-hand side of the equality.

Let us define the functions

$$K(q) = \int_0^{y_1} \frac{\tilde{y}_1}{D(\tilde{q}(\tilde{y}_1, t))} d\tilde{y}_1$$

$$K(p) = \int_0^{y_1} \frac{\tilde{y}_1}{D(\tilde{p}(\tilde{y}_1, t))} d\tilde{y}_1.$$

Then Lemma 3 can be written like

$$(7) \qquad \frac{d\delta}{dt} K(q) = S_2 - E(x_1, x_2, t),$$

and from Lemma 2

$$(8) \qquad S_1 = \left(\frac{d\delta}{dt} + \frac{d\beta}{dt}\right) K(p) + \frac{d\delta}{dt}[K(q) - K(p)] + O(\gamma).$$

In (7) and (8), we know that $S_2$ and $E$ are bounded for all time and that

$$\frac{d\gamma}{dt} = \frac{d\delta}{dt} + \frac{d\beta}{dt}.$$

We also know that that there are two positive constants $M$ and $c$ such that $M \geq K(q) \geq c > 0$, $M \geq K(p) \geq c > 0$ and

$$K(p) - K(q) = O(\gamma).$$

By combining (7) and (8) we obtain

$$S_1 = \frac{d\gamma}{dt}K(p) + O(\gamma).$$

Now we use Lemma 1 to get the estimate

(9) $$|\frac{d\gamma}{dt}| \leq (\text{const})|\gamma \cdot \log \gamma|,$$

if $\gamma$ is less than a small constant.

The proof of Theorems 1 and 2 follows directly from integrating (9).

## 7. Stretching formula

In the previous section we proved that if a solution $\theta$ of equation (2) is constant along hyperbolas then the angle of the saddle $\gamma$ cannot close in finite time. The angle cannot close faster than a double exponential in time.

In this section we want to prove Corollaries 1 and 2. Furthermore, if we assume that $\int_0^t \Phi(s)ds$ is bounded for finite time then we can rule out a simple hyperbolic blow-up scenario by the necessary and sufficient condition (3) mentioned in Section 2.

We use the formula

(10) $$\left(\frac{\partial}{\partial t} + u \cdot \nabla\right)|\nabla^\perp \theta| = \alpha|\nabla^\perp \theta|,$$

where $\alpha$ is the level set stretching factor

$$\alpha = \frac{1}{2}(\nabla u + \nabla^t u)\xi \cdot \xi.$$

Here, $\xi$ is the direction field of $\nabla^\perp \theta$. For (10), see Constantin [3] and Constantin, Majda and Tabak [7]. Consequently $\alpha$ can be represented by

$$\alpha(x) = \text{P.V.} \int_{R^2} ((\frac{y}{|y|} \cdot \xi^\perp(x))(\xi(x+y) \cdot \xi(x)))|\nabla^\perp \theta|(x+y)\frac{dy}{|y|^2}.$$

In [7] they derive an upper bound for the magnitude of the stretching factor

$$|\alpha(x)| \leq C[G(x)|u(x)| + (\tilde\tau G(x) + 1)(G||\theta||_{L^\infty} + \tilde\tau^{-2}||\theta||_{L^2})].$$

C, $\tilde\tau$ are fixed constants and

$$G(x) = \sup_{|y|\leq\tilde\tau}|\nabla \xi(x+y)|.$$

Now we derive an estimate for the velocity $u(x)$:

$$u(x) = \int_{R^2} \frac{y^\perp \theta(x+y)}{|y|^3}dy.$$

We consider $\tau > 0$ and we define $\aleph(r)$ to be a smooth nonnegative function such that $\aleph(r) = 1$ for $0 \leq r \leq 1$, and $\aleph(r) = 0$ for $r \geq 2$. Then we split $u(x)$ in two integrals:

$$u(x) = I_1 + I_2 = \int \aleph(\frac{|y|}{\tau}) \frac{y^\perp \theta(x+y)}{|y|^3} dy + \int (1 - \aleph(\frac{|y|}{\tau})) \frac{y^\perp \theta(x+y)}{|y|^3} dy.$$

By integrating by parts, we can easily prove that

$$|I_1| \leq C\tau \cdot \sup |\nabla^\perp \theta| + C.$$

For $I_2$ we obtain

$$|I_2| \leq \int_{|y|>2\tau} \frac{|\theta(x+y)|}{|y|^2} dy$$
$$= \int_{2\tau<|y|<k} + \int_{|y|\geq k}$$

where $k$ is fixed.

By performing the first integral and for the second one using the conservation of $\int_{R^2} |\theta|^2 dx$,

$$|I_2| \leq C|\log \tau| + \text{const}.$$

Therefore

$$|u(x)| \leq C_1 \sup |\nabla^\perp \theta| \tau + C_2 |\log \frac{1}{\tau}| + C_3.$$

Take $\tau = \frac{1}{\sup |\nabla^\perp \theta|}$. Then $u$ is bounded by $\log \sup |\nabla^\perp \theta|$, assuming that $\sup |\nabla^\perp \theta| > e$.

Under the assumptions of Theorem 1, $G(x)$ is bounded by a double exponential for $|x| > c$. Therefore we can estimate $\alpha$ on $V \equiv U \setminus \{|x| < c\}$ by

(11) $$|\alpha| \leq \log ||\nabla^\perp \theta||_{L^\infty} e^{e^t}.$$

PROPOSITION 1. *Assume the existence of a function $f = f(\rho, t)$ independent of $\sigma$. Then the material derivative ($\equiv D_t = \frac{\partial}{\partial t} + u \cdot \nabla$) of $f$ is also independent of $\sigma$.*

To prove the proposition we have to compute the material derivative:

$$D_t f = \left(\frac{\partial}{\partial t} + u \cdot \nabla\right) f(\rho, t) = \frac{\partial f}{\partial t} + \frac{\partial f}{\partial \rho}\left(\frac{\partial \rho}{\partial t} - \frac{\partial \psi}{\partial \sigma}\right).$$

From equation (5) and the fact that the functions $\frac{\partial f}{\partial t}$, $\frac{\partial f}{\partial \rho}$ are independent of $\sigma$, we can conclude that the material derivative of $f$ is also independent of $\sigma$.

In particular the function $F = \frac{|\nabla^\perp \theta|}{|\nabla \rho|}$ is independent of $\sigma$. We write the material derivative of $F$ as follows:

$$(12) \quad D_t\left(\frac{|\nabla^\perp \theta|}{|\nabla \rho|}\right) = \frac{1}{|\nabla \rho|} D_t(|\nabla^\perp \theta|) + |\nabla^\perp \theta| D_t\left(\frac{1}{|\nabla \rho|}\right)$$

$$= \left(\alpha + |\nabla \rho| D_t\left(\frac{1}{|\nabla \rho|}\right)\right) \frac{|\nabla^\perp \theta|}{|\nabla \rho|}.$$

We can estimate $|\nabla \rho| D_t(\frac{1}{|\nabla \rho|})$ on $V$ by

$$(13) \quad |\nabla \rho| |D_t\left(\frac{1}{|\nabla \rho|}\right)| = \frac{1}{|\nabla \rho|} |D_t(|\nabla \rho|)|$$

$$= \frac{1}{|\nabla \rho|} \left|\left(\frac{\partial}{\partial t} + u \cdot \nabla\right) |\nabla \rho|\right|$$

$$\leq C_1 \left|\frac{d\delta}{dt}\right| + C_2 \left|\frac{d\beta}{dt}\right| + C \log ||\nabla^\perp \theta||_{L^\infty} e^{e^t}.$$

Lemma 3 shows that $|\frac{d\delta}{dt}| \leq$ const. The same estimate can be obtained, $|\frac{d\beta}{dt}| \leq$ const, by substituting the point $q$ with $p$ in Lemma 3.

Combining (11) and (13) we obtain

$$\left|\alpha + |\nabla \rho| D_t\left(\frac{1}{|\nabla \rho|}\right)\right| \leq e^{e^t} \log ||\nabla^\perp \theta||_{L^\infty} + \text{const}$$

on $V$ when $||\nabla^\perp \theta||_{L^\infty} > e^{e^t}$. Then by (12)

$$(14) \quad |D_t F| \leq e^{e^t} \log ||\nabla^\perp \theta||_{L^\infty} \cdot F$$

on $V$ when $||\nabla^\perp \theta||_{L^\infty} > e^{e^t}$. However, because $F$ only depends on $\rho$ and $t$ and by the proposition it follows that inequality (14) only depends on $\rho$ and $t$.

If $|x| < c$, then find $\tilde{x} \in V$ with $\rho(x) = \rho(\tilde{x})$. Inequality (14) holds for $\tilde{x}$ and therefore it holds for $x$. Thus, (14) holds, not just on $V$, but on $U$.

We want to estimate the material derivative of $|\nabla^\perp \theta|$ on $U$. By (12) we have the following equality

$$(15) \quad D_t(|\nabla^\perp \theta|) = |\nabla \rho| \cdot D_t F + F \cdot D_t(|\nabla \rho|).$$

The estimate for the first term of the right-hand side of equality (15) is obtained in (14). For the second term we know that on $V$,

$$(16) \quad F \leq ||\nabla^\perp \theta||_{L^\infty}.$$

Because $F$ is independent of $\sigma$, inequality (16) holds on $U$. We compute the material derivative of $|\nabla \rho|$ and estimate it on $U$ as follows:

$$|D_t(|\nabla \rho|)| = \left|\frac{d}{dt}|\nabla \rho| + u \cdot \nabla(|\nabla \rho|)\right|$$

$$\leq C_1 \left|\frac{d\delta}{dt}\right| + C_2 \left|\frac{d\beta}{dt}\right| + |u \cdot \nabla(|\nabla \rho|)|$$

$$\leq C_1 \left|\frac{d\delta}{dt}\right| + C_2 \left|\frac{d\beta}{dt}\right| + \log ||\nabla^\perp \theta||_{L^\infty} e^{e^t}.$$

We recall again Lemma 3 and get the following estimate on $U$:

$$|D_t(|\nabla^\perp\theta|)| \leq C||\nabla^\perp\theta||_{L^\infty}\log||\nabla^\perp\theta||_{L^\infty}e^{e^t}.$$

If $x \in R^2 \setminus U$ then $|\nabla\xi| \leq \Phi(t)$. By using the upper bound of $\alpha$ we obtain

$$|D_t(|\nabla^\perp\theta|)| \leq C||\nabla^\perp\theta||_{L^\infty}\log||\nabla^\perp\theta||_{L^\infty}(\Phi(t) + e^{e^t})$$

on $x \in R^2 \setminus U$.

Therefore

$$(17) \qquad |\frac{d}{dt}||\nabla^\perp\theta||_{L^\infty}| \leq C||\nabla^\perp\theta||_{L^\infty}\log||\nabla^\perp\theta||_{L^\infty}(e^{e^t} + \Phi(t))$$

whenever $||\nabla^\perp\theta||_{L^\infty} > e^{e^t}$.

At time $t = 0$, $||\nabla^\perp\theta||_{L^\infty} < \infty$. We finish the proofs of Corollaries 1 and 2 by integrating (17).

*Remark* 5. In [4], Constantin shows that under the hypothesis of a "proper nondegenerate-self-similar Ansatz" the $||\nabla^\perp\theta||_{L^\infty}$ grows at most as a double exponential.

*Remark* 6. Suppose $\theta = \theta(\Pi(x,t),t)$ is constant along ellipses. We allow these ellipses to close according to

$$\Pi(x,t) = a(t) \cdot x_1^2 + b(t) \cdot x_2^2.$$

Whenever this happens in the numerical simulations, the norm of the gradient of $\theta$ does not grow fast. We can show that in this case the $|\nabla^\perp\theta|$ is bounded by a double exponential in time.

## 8. Similar results for Euler

Two dimensional case: Majda and Tabak [12] ran numerical simulations with the same initial data for the 2D Euler vorticity equation and compared them with the simulations in [7]. The norm of the gradient perpendicular of the vorticity grew much slower for 2D Euler than the $|\nabla^\perp\theta|$ above.

The representation of 2D Euler equation in vorticity form is

$$(\partial_t + u \cdot \nabla)\omega = 0$$
$$u = \nabla^\perp\psi \quad where \quad \omega = \triangle\psi.$$

The two active scalars $\theta$ and $\omega$ are similar, but they differ in the characterization of the stream function.

Using the same scheme as in Section 4, we assume that $\omega$ is constant along hyperbolas and that $\gamma$ is the angle of the saddle. We can show

$$|\log\gamma(t)| \leq (\text{constant}) \cdot t,$$

which means $\gamma$ can go to zero at most as an exponential. The proof is identical to the one in Section 6, but in Lemma 1 $\psi$ is defined by

$$\psi = \frac{1}{2\pi} \int_{R^2} \omega(x+y) \log|y| dy$$

and it follows that

$$|\psi(p) - \psi(q)| \leq K|\gamma|.$$

Three-dimensional case: As we explained in Section 2, the quasi-geostrophic equation is a two-dimensional model for the 3D incompressible Euler equation. The techniques used in this paper give analogous results for 3D Euler [10]:

THEOREM 3. *Let $u(x,t)$ be a smooth solution of the* 3D *Euler incompressible equation defined for $0 \leq t < T_*$, $x \in R^3$ with*

$$\omega = \mathrm{curl}(u) = \frac{|\omega|}{r}(-\frac{\partial \rho}{\partial x_2}, \frac{\partial \rho}{\partial x_1}, 0).$$

*Here, $r = |\nabla \rho|$, $u$ is bounded up to $t = T_*$ and $\rho$ is defined as in Theorem 1 with the same nonlinear time dependent coordinate change and the same assumptions. If $\tilde{\omega}$ is not zero in a disc in $U$, then $\lim_{t \to T_*} \gamma(t)$ exists, and is not $0$.*

COROLLARY 3. *Let $u$ be as in Theorem 3. Let $\xi = \frac{\omega}{|\omega|}$. Assume $|\nabla \xi| < C$ on $(R^3 \setminus U) \times [0, T_*]$. Then $u$ continues to some solution of the* 3D *Euler incompressible equation on $R^3 \times [0, T_* + \varepsilon]$ for some $\varepsilon > 0$.*

A similar estimate is obtained for the angle of the saddle $\gamma$. Note that vortex stretching may take place in the setting of Theorem 3.

PRINCETON UNIVERSITY, PRINCETON, NJ
*Current address*: THE INSTITUTE FOR ADVANCED STUDY, PRINCETON, NJ
*E-mail address*: dcg@ias.edu